# Characteristic $p$ Approaches to the Jacobian Conjecture

Jeffrey Lang, Mathematics Department, University of Kansas, Lawrence, KS 66045

ABSTRACT

We present several versions of the Jacobian Conjecture in characteristic $p > 0$ each of which if true would imply the Jacobian Conjecture in characteristic 0. We test these characteristic $p$ versions of the conjecture against several families of Jacobian pairs in characteristic $p$. Based on the results we propose a characteristic $p$ approach to solving the Jacobian Conjecture in characteristic 0.

## Section 1. Characteristic $p > 0$ Approaches to the Two-Dimensional Jacobian Conjecture

In this paper we present several versions of the Jacobian Conjecture in characteristic $p > 0$ each of which if true would imply the Jacobian Conjecture in characteristic 0. We then test these characteristic $p$ versions of the conjecture against several families of Jacobian pairs in characteristic $p$. When I began this project, I was very confident that among these families of Jacobian pairs I would ultimately find counterexamples to each of these characteristic $p$ versions of the Jacobian conjecture. To my surprise, I found that each of the characteristic $p$ versions of the Jacobian Conjecture hold for each of the families of Jacobians pairs considered in this article. Over the years I had gone back and forth on the efficacy of studying Jacobian pairs in characteristic $p > 0$ as regards the Jacobian Conjecture in characteristic 0. My last published views on the matter were on the pessimistic side.[1] The findings herein have me seriously questioning that judgement.

Throughout this article we let $k$ be a field and $k^*$ be the nonzero elements of $k$. For each positive integer $n$, let $A_n = k[x_1, \cdots, x_n]$ be the polynomial ring in $n$ variables over $k$. For $f_1, \cdots, f_m \in A_n$, let $k[f_1, \cdots, f_m]$ be the $k$−subalgebra of $k[x_1, \cdots, x_n]$ generated by the $f_i$ and $k(f_1, \cdots, f_m)$ be its quotient field. For $f_1, \cdots, f_m \in A_n$, let $J(f_1, \cdots, f_n)$ be the determinant of the $n \times n$ Jacobian matrix, $[\partial(f_i)/\partial x_j]$, $1 \leq i, j \leq n$. We say $(f_1, \cdots, f_n)$ is an **automorphic $n$−tuple** if $A_n = k[f_1, \cdots, f_n]$. Note that $(f_1, \cdots, f_n)$ is an automorphic $n$−tuple if and only if the $k$−endomorphism $\phi$ of $k[x_1, \cdots, x_n]$ defined by $\phi(x_i) = f_i, 1 \leq i \leq n$, is a $k$−algebra automorphism. If $(f_1, \cdots, f_n)$ is an automorphic $n$−tuple, then the Multivariable Chain Rule applied to the composition $\phi \circ \phi^{-1}$ implies $J(f_1, \cdots, f_n) \in k^*$. For $f_1, \cdots, f_n \in A_n$, we say $(f_1, \cdots, f_n)$ is a **Jacobian $n$−tuple** if $J(f_1, \cdots, f_n) \in k^*$. Thus, every automorphic $n$−tuple is a Jacobian $n$−tuple. The **Jacobian Conjecture**, introduced by O. H. Keller[2] in 1939, states the converse. If the characteristic of $k$ is 0 and $n = 1$, the Jacobian Conjecture is true. It is still unknown if it is true when the characteristic of $k$ is 0 and $n > 1$. The Jacobian Conjecture is not true when the characteristic of $k$ equals $p > 0$, even when $n = 1$. For example, when $p > 0$, the derivative of $f_1 = x_1^p + x_1$ equals 1, but $k[x_1] \neq k[f_1]$ since $[k(x_1):k(f_1)] = p$.

In this paper we will focus on the two variable Jacobian Conjecture, i.e. the case when $n = 2$. In this case we will refer to automorphic 2-tuples and Jacobian 2-tuples as **automorphic pairs** and **Jacobian pairs**, respectively. Since the Jacobian Conjecture is not true when the characteristic of $k$ is greater than 0, for this and other reasons many researchers of the $2$−Dimensional Jacobian Conjecture often bypass consideration of Jacobian pairs in positive characteristic. Yet their study might reveal important information about Jacobian pairs in characteristic 0. To start with, we have the following result relating Jacobian pairs in characteristic 0 to those in characteristic $p > 0$.

**Proposition 1.2.** *Let $f_1, f_2 \in \mathbb{C}[x_1, x_2]$. If $f_1, f_2$ is a Jacobian pair, then for all but a finite number of prime numbers $p > 0$, there exists a finite field $F$ of characteristic $p$ and a Jacobian pair $\widehat{f_1}, \widehat{f_2} \in F[x_1, x_2]$ such that the support of $f_i$ equals the support of $\widehat{f_i}$ for each $i = 1, 2$. In particular, the Newton polygons of $f_i$ and $\widehat{f_i}$ will be identical for each $i$.*[3](Proposition (1.2))

We also have the following equivalent formulations of the two variable Jacobian Conjecture in characteristic 0.

**Theorem 1.3.** (Abhyankar). *Let $k$ be a field of characteristic 0. Then the following statements are equivalent.*

(i) *If $f_1, f_2 \in k[x_1, x_2]$ and $\partial(f_1, f_2)/\partial(x_1, x_2) \in k^*$, then $k[f_1, f_2] = k[x_1, x_2]$.*
(ii) *If $f_1, f_2 \in k[x_1, x_2]$ and $\partial(f_1, f_2)/\partial(x_1, x_2) \in k^*$, then $f_i$ has one point at infinity (See Definition 2.15 below) for each $i$.*
(iii) *If $f_1, f_2 \in k[x_1, x_2]$ and $\partial(f_1, f_2)/\partial(x_1, x_2) \in k^*$, then the Newton-Polygon of $f_i$ is a triangle with vertices $(n, 0), (0, m)$ and $(0, 0)$ for some nonnegative integers $n$ and $m$.*
(iv) *If $f_1, f_2 \in k[x_1, x_2]$ and $\partial(f_1, f_2)/\partial(x_1, x_2) \in k^*$, then $\deg(f)$ divides $\deg(g)$ or $\deg(g)$ divides $\deg(f)$.*[4](Theorem 19.4)

Note that the conclusions of each of the four equivalent statements of Theorem 1.3 are properties of automorphic pairs in every characteristic. It follows from the above two results that each of the following conjectures, if true, would imply the 2-dimensional Jacobian Conjecture in characteristic 0.

**Low Degree Jacobian Conjecture 1.** Let $k$ be of characteristic $p > 0$ and $f_1, f_2 \in k[x_1, x_2]$ with $\deg(f_i) < p$ for each $i$. Then $\partial(f_1, f_2)/\partial(x_1, x_2) \in k^*$ if and only if $k[x_1, x_2] = k[f_1, f_2]$.

**Low Degree Jacobian Conjecture 2.** Let $k$ be of characteristic $p > 0$ and $f_1, f_2 \in k[x_1, x_2]$ with $\deg(f_i) < p$ for each $i$. Then $\partial(f_1, f_2)/\partial(x_1, x_2) \in k^*$ if and only if $f_1$ and $f_2$ have one point at infinity.

Likewise, we have two other low degree versions of the two variable Jacobian Conjecture corresponding to statements (iii) and (iv) of Theorem 1.3 which, if true, would imply the 2-dimensional Jacobian Conjecture in characteristic 0.

We can obtain yet another characteristic $p > 0$ approach to the 2-dimensional Jacobian Conjecture in characteristic 0 based on the Separable Jacobian Conjecture introduced by Kossivi Adjamagboo. He argues for extending the $n$-Dimensional Jacobian Conjecture to characteristic $p > 0$ by adding the hypothesis that $p$ does not divide the degree of the field extension $k(x_1, \cdots, x_n)/k(f_1, \cdots, f_n)$.[5] This translates in the two-variable case into the following.

**Two-Dimensional Separable Jacobian Conjecture in Characteristic $p > 0$.** Let $k$ be of characteristic $p > 0$ and $f_1, f_2 \in k[x_1, x_2]$. Then $k[x_1, x_2] = k[f_1, f_2]$ if and only if $\partial(f_1, f_2)/\partial(x_1, x_2) \in k^*$ and $[k(x_1, x_2): k(f_1, f_2)]$ is not divisible by $p$.

To show that the statement of this conjecture implies the 2-dimensional Jacobian Conjecture in characteristic 0 will require some preliminary results.

**Proposition 1.4.** *If $f_1, f_2 \in k[x_1, x_2]$ and $\partial(f_1, f_2)/\partial(x_1, x_2) \neq 0$, then $k(x_1, x_2)$ is separable algebraic over $k(f_1, f_2)$.*

**Proof.** It is enough to show that if $D$ is a derivation of $k(x_1, x_2)$ such that $k(f_1, f_2) \subset D^{-1}(0)$, then $D = 0$. [6](Proposition 8.18) If $D$ is such a derivation, then for each $i = 1, 2$, $D(f_i) = (\partial(f_i)/\partial x_1)D(x_1) + (\partial(f_i)/\partial x_2)D(x_2) = 0$. Since $\partial(f_1, f_2)/\partial(x_1, x_2) \neq 0$, $D(x_1) = D(x_2) = 0$. Hence, $D = 0$.

**Proposition 1.5.** *If $f_1, f_2 \in k[x_1, x_2]$ and $\partial(f_1, f_2)/\partial(x_1, x_2) \in k^*$, then $[k(x_1, x_2): k(f_1, f_2)] \leq \deg(f_1) \cdot \deg(f_2)$.*

**Proof.** By Proposition 1.4, $k(x_1, x_2)$ is a finite separable extension of $k(f_1, f_2)$. By the Primitive Element Theorem, there exists $\alpha \in k$ such that $k(x_1, x_2) = k(f_1, f_2, x_1 + \alpha x_2)$. After a homogeneous change of variables, we may assume that $k(x_1, x_2) = k(f_1, f_2, x_1)$.
 Let $\{t_i, t_2\}$ be a set of algebraic indeterminates over $k(x_1, x_2)$. For $i = 1, 2$, let $F_i(t_1, t_2) = f_i(t_1, t_2) - f_i(x_1, x_2)$. Then $F_1, F_2 \in k(f_1, f_2)[t_1, t_2]$ and are relatively prime since $\partial(F_1, F_2)/\partial(t_1, t_2) \in k^*$. Let $R(t_1) = \text{res}_{t_2}(F_1, F_2)$ (i.e. the resultant with respect to $t_2$). Then $R(t_1) \in k(f_1, f_2)[t_1]$ and $R(x_1) = 0$, since $F_i(x_1, x_2) = 0$ for each $i$. If we let $M(t_1) \in k(f_1, f_2)[t_1]$ be the minimal polynomial of $x_1$ over $k(f_1, f_2)$, it follows that $M(t_1)$ is a factor of $R(t_1)$. Hence, $\deg_{t_1}(M(t_1)) \leq \deg_{t_1}(R(t_1))$. Since $\deg_{t_1}(M(t_1)) = [k(x_1, x_2): k(f_1, f_2)]$ and $\deg_{t_1}(R(t_1)) \leq \deg(f_1)\deg(f_2)$, we have, $[k(x_1, x_2): k(f_1, f_2)] \leq \deg(f_1)\deg(f_2)$.

**Proposition 1.6.** *If the 2-Dimensional Separable Jacobian Conjecture in Characteristic $p > 0$ is true, then the 2-Dimensional Jacobian Conjecture in Characteristic 0 is as well.*

**Proof.** Assume that the 2-Dimensional Separable Jacobian Conjecture holds and let $f_1, f_2 \in \mathbb{C}[x_1, x_2]$ be a Jacobian pair. Then for all but a finite number of prime numbers $p > 0$, there exists a finite field $F$ of characteristic $p$ and a Jacobian pair $\widehat{f_1}, \widehat{f_2} \in F[x_1, x_2]$ such that the support of $f_i$ equals the support of $\widehat{f_i}$ for each i $= 1, 2$. Then for those $p$ we have $\deg(f_i) = \deg(\widehat{f_i})$ and for those $p$ that are also greater than $\deg(f_1) \cdot \deg(f_2)$, we have $[F(x_1, x_2): F(\widehat{f_1}, \widehat{f_2})]$ is not divisible by $p$ by Proposition 1.5. Since the Separable Jacobian Conjecture holds there exists a positive integer $N$ such that for all prime numbers $p > N$, there exists a finite field $F$ of characteristic $p$ and an automorphic pair $\widehat{f_1}, \widehat{f_2} \in F[x_1, x_2]$ such that the support of $f_i$ equals the support of $\widehat{f_i}$ for each i $= 1, 2$, which implies $\deg(f_1)$ divides $\deg(f_2)$ or vice versa. It follows from Theorem 1.3 that every Jacobian pair in $\mathbb{C}[x_1, x_2]$ is an automorphic pair.

Another feature of Jacobian pairs in characteristic $p > 0$ that might prove useful is that we know a fair amount about them structurally. Let $\nabla$ be the differential operator on $k(x_1, \cdots, x_n)$ defined by $\nabla = (\partial^{n(p-1)})/(\partial x_1^{p-1} \cdots \partial x_n^{p-1})$. For each $i = 1, 2, \cdots, n$, let $D_i$ be the $k$-derivation on $k(x_1, \cdots, x_n)$ defined by $D_i(g) = J(f_1, \cdots, f_{i-1}, g, f_{i+1}, \cdots, f_n)$ for all $g \in k(x_1, \cdots, x_n)$. The next two theorems give equivalent characterizations of Jacobian $n$-tuples in characteristic $p > 0$.

**Theorem 1.7.** (Nousainen) *Let $f_1, \cdots, f_n \in A_n$. Then the following conditions are equivalent.*

(1) $(f_1, \cdots, f_n)$ *is a Jacobian $n$-tuple.*
(2) $k[x_1, \cdots, x_n] = k[x_1^p, \cdots, x_n^p, f_1, \cdots, f_n]$.
(3) *The monomials, $f_1^{i_1} \cdots f_n^{i_n}$, $0 \le i_j \le p-1$, form a free basis of the $k[x_1^p, \cdots, x_n^p]$-module $k[x_1, \cdots, x_n]$.*[7](See Theorem 2.1)

**Theorem 1.8.** (Lang and Mandal) *Let $f_1, \cdots, f_n \in k[x_1, \cdots, x_n]$. Then the following conditions are equivalent.*

1. $(f_1, \cdots, f_n)$ is a Jacobian $n$-tuple.
2. There exists $\alpha \in k^*$ such that for each $i = 1, 2, \cdots, n$, and each $g \in k(x_1, \cdots, x_n)$, $g = \alpha \sum_{j=0}^{p-1} f_i^j D_i^{p-1}(f_i^{p-j-1} g)$.
3. $\nabla = \alpha D_1^{p-1} \cdots D_n^{p-1}$ for some $\alpha \in k^*$.
4. There exists $\alpha \in k^*$ such that $\nabla(f_1^{r_1} \cdots f_n^{r_n}) = \begin{cases} 0, & \text{if } 0 \le r_i < p - 1, \text{ for some } i. \\ \alpha, & \text{if each } r_i = p - 1. \end{cases}$

[8] (Theorem 2.2)

We also have the following identity that explicitly relates the Jacobian determinant, $J(f_1, \cdots, f_n)$, for $f_1, \cdots, f_n \in A_n$, to a specific element of $k[x_1^p, \cdots, x_n^p, f_1, \cdots, f_n]$.

**Theorem 1.9.** (Lang) *If the characteristic of $k$ is $p > 0$ and $f_1, \cdots, f_n \in k[x_1, \cdots, x_n]$, then*

$$\sum_{i_1=0}^{p-1} \cdots \sum_{i_n=0}^{p-1} f_1^{i_1} \cdots f_n^{i_n} \nabla(f_1^{p-1-i_1} \cdots f_n^{p-1-i_n}) = (-1)^n (J(f_1, \cdots, f_n))^{p-1}.$$ [9](Theorem 1.7)

Based on the above considerations we would like to propose the following admittedly characteristic $p > 0$ approach to the Two-dimensional Jacobian Conjecture.

Attempt to describe all methods of generating Jacobian pairs in characteristic $p > 0$. If no such pair serves as a counterexample to the Low Degree Jacobian Conjecture 1, then the Two-dimensional Jacobian Conjecture in characteristic 0 is true. Likewise, if no such pair serves as counterexample to the Low Degree Jacobian Conjecture 2 or to the Separable Jacobian Conjecture then the Two-dimensional Jacobian Conjecture in characteristic 0 is true. Of course, describing all methods for generating Jacobian pairs in positive characteristic could be a formidable task, but our hope is, based partly on the preceding three theorems, that their number is relatively small and that this goal is within reach. Be that as it may, it is still a worthwhile problem that whenever we do produce a genuinely new algorithm for generating Jacobian pairs, we should test the pairs it produces against the above Jacobian conjectures in positive characteristic. This we do below for several such algorithms known to us.

**Section 2. The Separable Jacobian Conjecture and $p-$Morphisms.**

In 1942 Jung[10](section 4) first proved that every automorphism of the polynomial ring in two variables over a field of characteristic 0 is a composition of elementary automorphisms. In 1953 Van der Kulk[11] extended Jung's theorem to a field of arbitrary characteristic. Since then several other proofs have appeared.

**Theorem 2.1. Automorphism Theorem of $k[x_1, x_2]$.** *Every $k-$ algebra automorphism of $k[x_1, x_2]$ is a finite compositional product of automorphisms of the type: (1) $x_1 \to a_{11}x_1 + a_{12}x_2 + a_{13}, x_2 \to a_{21}x_1 + a_{22}x_2 + a_{23}$ with the $a_{ij} \in k$ and $a_{11}a_{22} - a_{12}a_{21} \neq 0$ and (2) $x_1 \to x_1, x_2 \to x_2 + h(x_1)$ with $h(x_1) \in k[x_1]$.*

Thus, if the Two-Dimensional Jacobian Conjecture is true in characteristic 0 then every Jacobian pair will be of the form $(\phi(x_1), \phi(x_2))$ where $\phi$ is a product of elementary automorphisms, i.e. type 1 and type 2 automorphisms as described in Theorem 2.1. Even if the Two-Dimensional Jacobian Conjecture is not true in characteristic 0, all finite products of elementary automorphisms will produce Jacobian pairs in this way.
      The natural extension of this method of generating Jacobian pairs in characteristic 0 to characteristic $p > 0$ is to replace $h(x_1)$ in the definition of type 2 automorphism with an element in the kernel of $\partial/\partial x_2$, which will be of the form $h(x_1, x_2^p)$ for some $h(x_1, x_2) \in k[x_1, x_2]$. We will call this type of $k-$endomorphism of $k[x_1, x_2]$ a **type 2\*** endomorphism. A finite product of type 1 and type 2* maps will in general not be an automorphism of $k[x_1, x_2]$ but it will be a $k-$endomorphism $\phi$ of $k[x_1, x_2]$ such that $(\phi(x_1), \phi(x_2))$ is a Jacobian pair.
      From here on in this paper we will assume, unless stated otherwise, that the characteristic of $k$ is $p > 0$.

**Definition 2.2.** A map $\phi: k[x_1, x_2] \to k[x_1, x_2]$ is a $p-$ **morphism** of $k[x_1, x_2]$ if it is a finite compositional product of $k-$endomorphisms of the type: (1) $x_1 \to a_{11}x_1 + a_{12}x_2 + a_{13}, x_2 \to a_{21}x_1 + a_{22}x_2 + a_{23}$ with the $a_{ij} \in k$ and $a_{11}a_{22} - a_{12}a_{21} \neq 0$ and (2*) $x_1 \to x_1, x_2 \to x_2 + h(x_1, x_2^p)$ with $h(x_1, x_2) \in k[x_1, x_2]$. By the Multivariable Chain Rule, if $\phi$ is a $p-$ morphism of $k[x_1, x_2]$, then $(\phi(x_1), \phi(x_2))$ will be a Jacobian pair.

**Remark 2.3.** If $f_1, f_2 \in k[x_1, x_2]$ and $\partial(f_1, f_2)/\partial(x_1, x_2) \neq 0$, then $f_1$ and $f_2$ are algebraically independent over $k$. Thus, in this case, $k[f_1, f_2]$ will be isomorphic to a polynomial ring over $k$ in two variables. Since $k(f_1, f_2) \subset k(x_1, x_2)$ and both fields have transcendence degrees two over $k$, $k(x_1, x_2)$ will be algebraic over $k(f_1, f_2)$. In the case where $\partial(f_1, f_2)/\partial(x_1, x_2) \in k^*$, we will also have that the greatest common divisor of $f_1$ and $f_2$ in $k[x_1, x_2]$ will be 1. For if $g$ is a common factor of $f_1$ and $f_2$ in $k[x_1, x_2]$, then $g$ would be a factor of $\partial(f_1, f_2)/\partial(x_1, x_2)$; i.e. $g$ would belong to $k^*$.

**Proposition 2.4.** *Let $f_1, f_2, u_1, u_2 \in k[x_1, x_2]$. If the pairs $f_1, f_2$ and $u_1, u_2$ are algebraically independent over $k$, then $k(x_1, x_2)$ is algebraic over $k(f_1(u_1, u_2), f_2(u_1, u_2))$ and*

$$[k(x_1, x_2): k(f_1(u_1, u_2), f_2(u_1, u_2))]$$
$$= [k(x_1, x_2): k(u_1, u_2)][k(x_1, x_2): k(f_1(x_1, x_2), f_2(x_1, x_2))].$$

**Proof.** By Remark 2.5, $k(x_1, x_2)$ is algebraic over $k(u_1, u_2)$ and $k(u_1, u_2)$ is algebraic over $k(f_1(u_1, u_2), f_2(u_1, u_2))$. Hence, $k(x_1, x_2)$ is algebraic over $k(f_1(u_1, u_2), f_2(u_1, u_2))$ and

$$\big[[k(x_1, x_2): k(f_1(u_1, u_2), f_2(u_1, u_2))]\big] =$$
$$[k(x_1, x_2): k(u_1, u_2)][k(u_1, u_2): k(f_1(u_1, u_2), f_2(u_1, u_2))]$$

We also have the following commutative diagram of maps, where the upper horizontal map is the $k$−isomorphism that sends $x_i$ to $u_i$ for each $i$ and the lower horizontal map is the isomorphism induced by the upper horizontal map.

$$\begin{array}{ccc} k(x_1, x_2) & \to & k(u_1, u_2) \\ \cup & & \cup \\ k(f_1(x_1, x_2), f_2(x_1, x_2)) & \to & k(f_1(u_1, u_2), f_2(u_1, u_2)) \end{array}$$

It follows that $[k(u_1, u_2): k(f_1(u_1, u_2), f_2(u_1, u_2))] = [k(x_1, x_2): k(f_1(x_1, x_2), f_2(x_1, x_2))]$, from which the statement of the Proposition follows.

**Definition 2.5.** If $\phi$ is a $k$−endomorphism of $k[x_1, x_2]$ such that $\phi(x_1)$ and $\phi(x_2)$ are algebraically independent over $k$, we let $\deg(\phi) = [k(x_1, x_2): k(\phi(x_1), \phi(x_2))]$.

**Corollary 2.6.** *For each $i = 1, 2, \cdots, n$, assume $\phi_i$ is a $k$−endomorphism of $k[x_1, x_2]$ such that $\phi_i(x_1)$ and $\phi_i(x_2)$ are algebraically independent over $k$ and let $\phi = \phi_1 \phi_2 \cdots \phi_n$. Then $\deg(\phi) = \prod_{i=1}^{n} \deg(\phi_i)$.*

**Proof.** The $n = 2$ case is a direct consequence of Proposition 2.4. The general case follows by induction on $n$.

**Definition 2.7.** For $\phi$ a type $2^*$ morphism of $k[x_1, x_2]$ defined by $\phi(x_1) = x_1$ and $\phi(x_2) = x_2 + h(x_1, x_2^p)$, with $h(x_1, x_2) \in k[x_1, x_2]$, we define $\deg_p \phi = \deg_{x_2}(h(x_1, x_2))$. In the case where $\phi$ is a type 1 map, we define $\deg_p \phi = 0$.

**Remark 2.8.** Note that if $\phi$ is a type $2^*$ morphism of $k[x_1, x_2]$ and $\deg_p \phi = 0$, then $\phi$ is in fact a type 2 automorphism as described in Theorem 2.1.

**Proposition 2.9.** *Let $\phi$ be a $p$ − morphism of $k[x_1, x_2]$ with $\phi = \phi_1 \phi_2 \cdots \phi_n$ where the $\phi_i$ are type 1 and type $2^*$ morphisms as described in Definition 2.2. Then,*

$$\deg(\phi) = p^{\left(\sum_{i=1}^{n} \deg_p \phi_i\right)}.$$

**Proof.** By Corollary 2.6. it is enough to show that if $\phi$ is a type 1 or type $2^*$ morphism, then $\deg(\phi) = p^{\deg_p \phi}$. If $\phi$ is a type 1 morphism defined by $\phi(x_1) = a_{11} x_1 + a_{12} x_2 + a_{13}$ and $\phi(x_2) = a_{21} x_1 + a_{22} x_2 + a_{23}$ with the $a_{ij} \in k$ and $a_{11} a_{22} - a_{12} a_{21} \neq 0$, then $k(\phi(x_1), \phi(x_2)) = k[x_1, x_2]$. Hence, $\deg(\phi) = 1 = p^{\deg_p \phi}$, since $\deg_p \phi = 0$ by definition. If $\phi$ is a type $2^*$ morphism defined by $\phi(x_1) = x_1$ and $\phi(x_2) = x_2 + h(x_1, x_2^p)$ with $h(x_1, x_2) \in k[x_1, x_2]$, then $k(\phi(x_1), \phi(x_2)) = k(x_1, x_2 + h(x_1, x_2^p))$ and $k(x_1, x_2)$ is a primitive extension of $k(\phi(x_1), \phi(x_2))$ generated by $x_2$. Let $t$ be a variable over $k(x_1, x_2)$ and let $m(t) = t + h(\phi(x_1), t^p) - \phi(x_2)$ in $k(\phi(x_1), \phi(x_2))[t]$. Then $m(t) = t + h(x_1, t^p) - (x_2 + h(x_1, x_2^p))$ and $k[\phi(x_1), \phi(x_2), t]/(m(t))$ is isomorphic to $k[x_1, t]$. Hence, $m(t)$ is irreducible in $k[\phi(x_1), \phi(x_2), t]$ and, thus, irreducible over $k(\phi(x_1), \phi(x_2))$ by Gauss's Lemma. Since $x_2$ is a root of $m(t)$, we obtain $[k(x_1, x_2): k(\phi(x_1), \phi(x_2))] = \deg(m(t)) = p^{\deg_{x_2}(h(x_1, x_2))} = p^{\deg_p \phi}$ by definition. Hence, $\deg(\phi) = p^{\deg_p \phi}$.

**Corollary 2.10.** *Let $\phi$ be a $p-$morphism of $k[x_1, x_2]$. Then $\phi$ is an automorphism if and only if $\deg(\phi) \not\equiv 0 \pmod{p}$.*

**Proof.** We have for some positive integer $n$ that $\phi = \phi_1 \phi_2 \cdots \phi_n$ where the $\phi_i$ are type 1 and type $2^*$ morphisms as described in Definition 2.2. If $\phi$ is an automorphism, then $[k(x_1, x_2): k(\phi(x_1), \phi(x_2))] = 1$. Hence, $[k(x_1, x_2): k(\phi(x_1), \phi(x_2))] \not\equiv 0 \pmod{p}$. Conversely, if $[k(x_1, x_2): k(\phi(x_1), \phi(x_2))] \not\equiv 0 \pmod{p}$, then by Proposition 2.9, $\deg_p \phi_i = 0$, for each $i$, which implies each $\phi_i$ is a type 1 or 2 automorphism as described in Theorem 2.1 (see Remark 2.8). Hence, $\phi$ is an automorphism.

**Remark 2.11.** If $\phi$ is a $p-$morphism, then $(\phi(x_1), \phi(x_2))$ is a Jacobian pair. If, conversely, each Jacobian pair is $(\phi(x_1), \phi(x_2))$ for some $p-$morphism $\phi$, then the Two-Variable Separable Jacobian Conjecture in positive characteristic is true by Corollary 2.10 and then the two-variable Jacobian Conjecture in characteristic 0 is true by Proposition 1.6.

Since we will not find counterexamples to the Separable Jacobian Conjecture in positive characteristic among pairs of the form $(\phi(x_1), \phi(x_2))$ where $\phi$ is a $p-$morphism of $k[x_1, x_2]$, two questions arise. (1) Is every Jacobian pair in characteristic $p > 0$ equal to $(\phi(x_1), \phi(x_2))$ for some $p-$morphism $\phi$ of $k[x_1, x_2]$? (2) Do the Low Degree Jacobian Conjectures hold for Jacobian pairs produced by $p-$morphisms? That is, if $\phi$ is a $p-$morphism of $k[x_1, x_2]$ and $\deg(\phi(x_i)) < p$ for each $i$, must $(\phi(x_1), \phi(x_2))$ be an automorphic pair or must $\phi(x_1)$ and $\phi(x_2)$ have one point at infinity? It turns out that these two questions are closely related. We will start by considering the number of points at infinity of Jacobian pairs produced by $p-$morphisms.

**Lemma 2.12.** *Let $F_1, F_2 \in k[x_1, x_2]$ be homogeneous polynomials of positive degrees $d_1, d_2$, respectively, such that $\partial(F_1, F_2)/\partial(x_1, x_2) = 0$. If $d_1 d_2 \not\equiv 0 \pmod{p}$, then every linear factor of $F_1$ of multiplicity not divisible by $p$ is also a factor of $F_2$ of multiplicity not divisible by $p$. If $d_1 \equiv 0 \pmod{p}$, then $d_2 \equiv 0 \pmod{p}$ or $F_1 \in k[x_1^p, x_2^p]$.*

**Proof.** Since $k$ is algebraically closed, for each $i = 1, 2$, we may factor $F_i$ as $F_i(x_1, x_2) = G_i(x_1, x_2)\bigl(H_i(x_1, x_2)\bigr)^p$, where $G_i(x_1, x_2), H_i(x_1, x_2) \in k[x_1, x_2]$ are homogeneous polynomials such that $G_i(x_1, x_2)$ and $H_i(x_1, x_2)$ are relatively prime and each of the linear factors of $G_i$ have multiplicity not divisible by $p$. Since $\partial(F_1, F_2)/\partial(x_1, x_2) = 0$, we obtain by Euler's Homogeneous Function Theorem,

$$0 = x_1\bigl((F_1)_{x_1}(F_2)_{x_2} - (F_1)_{x_2}(F_2)_{x_1}\bigr) + x_2\bigl((F_1)_{x_2}(F_2)_{x_2} - (F_1)_{x_2}(F_2)_{x_2}\bigr) = d_1 F_1 (F_2)_{x_2} - d_2 (F_1)_{x_2} F_2,$$

and

$$0 = x_2\bigl((F_1)_{x_1}(F_2)_{x_2} - (F_1)_{x_2}(F_2)_{x_1}\bigr) + x_1\bigl((F_1)_{x_1}(F_2)_{x_1} - (F_1)_{x_1}(F_2)_{x_1}\bigr) = d_2 F_2 (F_1)_{x_1} - d_1 (F_2)_{x_1} F_1.$$

Hence, if $d_1 d_2 \not\equiv 0 \pmod{p}$, then $\partial/\partial x_1 ((F_2)^{d_1}/(F_1)^{d_2}) = \partial/\partial x_2 ((F_2)^{d_1}/(F_1)^{d_2}) = 0$. From this it follows that $(F_2)^{d_1}/(F_1)^{d_2} \in k(x_1^p, x_2^p)$. Thus, there exists relatively prime homogeneous polynomials, $Q_1, Q_2 \in k[x_1, x_2]$, of the same degree such that $\bigl(F_1(x_1, x_2)\bigr)^{d_2}\bigl(Q_1(x_1, x_2)\bigr)^p = \bigl(F_2(x_1, x_2)\bigr)^{d_1}\bigl(Q_2(x_1, x_2)\bigr)^p$. Hence,

$$\bigl(G_1(x_1,x_2)\bigr)^{d_2}\bigl(H_1(x_1,x_2)\bigr)^{pd_2}\bigl(Q_1(x_1,x_2)\bigr)^p = \bigl(G_2(x_1,x_2)\bigr)^{d_1}\bigl(H_2(x_1,x_2)\bigr)^{pd_1}\bigl(Q_2(x_1,x_2)\bigr)^p.$$

If $u$ is a linear homogeneous factor of $F_1$ of multiplicity $r \not\equiv 0 \pmod{p}$, then $u$ is a factor of $(G_1)^{d_2}$ of multiplicity $rd_2 \not\equiv 0 \pmod{p}$, since $d_2 \not\equiv 0 \pmod{p}$. Hence, $u$ is a factor of $(G_1)^{d_2}(H_1)^{pd_2}(Q_1)^p$ of multiplicity that is congruent to $rd_2 \not\equiv 0 \pmod{p}$, from which it follows that $u$ must be a factor of $G_2$ and, therefore, a factor of $F_2$ of multiplicity not divisible by $p$.

If $d_1 \equiv 0$, then $d_2(F_1)_{x_2} F_2 = d_2 F_2 (F_1)_{x_1} = 0$, since $d_1 F_1 (F_2)_{x_2} - d_2 (F_1)_{x_2} F_2 = d_2 F_2 (F_1)_{x_1} - d_1 (F_2)_{x_1} F_1 = 0$. Hence, either $d_2 \equiv 0 \pmod{p}$ or $(F_1)_{x_1} = (F_1)_{x_2} = 0$, which implies $F_1 \in k[x_1^p, x_2^p]$.

**Notation 2.13.** Let $\mathbb{N}_0$ denote the set of nonnegative integers. Let $f \in k[x_1, x_2]$ be nonzero. Then $f = \sum_{j \in \mathbb{N}_0} f_j$, where each $f_j \in k[x_1, x_2]$ is nonzero homogeneous of degree $j$ or $f_j = 0$. We let $f^+$ denote the nonzero homogeneous form of $f$ of highest degree. If $f \notin k[x_1^p, x_2^p]$, then there exists a positive integer $m$ such that $f_m \notin k[x_1^p, x_2^p]$ and $f_j \in k[x_1^p, x_2^p]$ for all $j > m$. In this case we let $\bar{f} = \sum_{j=0}^m f_j$. Then $\bar{f}^+ = f_m$, i.e. $\bar{f}^+$ equals the highest degree form of $f$ that does not belong to $k[x_1^p, x_2^p]$.

**Remark 2.14.** Let $f_1, f_2 \in k[x_1, x_2]$. Then $\partial(f_1, f_2)/\partial(x_1, x_2) = \partial(\bar{f}_1, \bar{f}_2)/\partial(x_1, x_2)$ since $f_i - \bar{f}_i \in k[x_1^p, x_2^p]$ for each $i$. Hence, $(f_1, f_2)$ is a Jacobian pair if and only if $(\bar{f}_1, \bar{f}_2)$ is a Jacobian pair.

**Definition 2.15.** Let $f$ be nonzero in $k[x_1, x_2]$. We say that $f$ **has $r$ points at infinity** if $f^+ = u_1^{e_1} u_2^{e_2} \cdots u_r^{e_r}$ where the $u_i \in k[x_1, x_2]$ are mutually coprime linear homogeneous and the $e_i$ are positive integers. We say $f$ has **$n$ points at infinity modulo $p$** if $\bar{f}^+ = u_1^{e_1} u_2^{e_2} \cdots u_r^{e_r} H$, where the $u_i \in k[x_1, x_2]$ are mutually coprime linear homogeneous, $H \in k[x_1^p, x_2^p]$ is homogeneous, and the $e_i$ are positive integers not divisible by $p$.

**Proposition 2.16.** *Let $f_1, f_2 \in k[x_1, x_2]$ have positive degrees $d_1, d_2$, respectively, with $d_1 + d_2 > 2$. Assume $\deg(\partial(f_1, f_2)/\partial(x_1, x_2)) < d_1 + d_2 - 2$. If $d_1 d_2 \not\equiv 0 \pmod{p}$, then $f_1$ and $f_2$ have the same number of points at infinity modulo $p$ associated with the same linear homogeneous factors in $k[x_1, x_2]$. If $d_1 \equiv 0 \pmod{p}$, then $d_2 \equiv 0 \pmod{p}$ or $f_1^+ \in k[x_1^p, x_2^p]$.*

**Proof.** We have $\deg(\partial(f_1, f_2)/\partial(x_1, x_2)) \leq d_1 + d_2 - 2$ and $\deg(\partial(f_1, f_2)/\partial(x_1, x_2)) < d_1 + d_2 - 2$ if and only if $\partial(f_1^+, f_2^+)/\partial(x_1, x_2) = 0$. The conclusion then follows by Lemma 2.12.

**Lemma 2.17.** *Let $f_1, f_2 \in k[x_1, x_2]$ and $\partial(f_1, f_2)/\partial(x_1, x_2) \in k^*$. If $f_1$ has one point at infinity modulo $p$, then $f_2$ has one point at infinity modulo $p$.*

Proof. We have $(\bar{f}_1, \bar{f}_2)$ is a Jacobian pair (Remark 2.14). If $\deg(\bar{f}_1) + \deg(\bar{f}_2) = 2$, then $\bar{f}_1$ and $\bar{f}_2$ are coprime and linear in $k[x_1, x_2]$ so that each $f_i$ has one point at infinity modulo $p$. So we may assume $\deg(\bar{f}_1) + \deg(\bar{f}_2) > 2$. We have $\deg(\bar{f}_1) \not\equiv 0 \pmod{p}$ since $f_1$ has one point at infinity modulo $p$. We also have $\bar{f}_2^+ \notin k[x_1^p, x_2^p]$ by definition, which implies $\deg(\bar{f}_2) \not\equiv 0 \pmod{p}$ by Proposition 2.16. Hence, $f_2$ has one point at infinity modulo $p$ by Proposition 2.16.

**Theorem 2.18.** *Let $\phi$ be a $p-$ morphism of $k[x_1, x_2]$. Then for each $i = 1, 2$, $\phi(x_i)$ has one point at infinity modulo $p$.*

**Proof**. By Lemma 2.17, it is enough to show that at least one of $\phi(x_1)$ and $\phi(x_2)$ has one point at infinity modulo $p$. We will do so by proceeding by induction on the number of type 1 and type $2^*$ maps in the compositional product that defines $\phi$.

If $\phi$ is a type 1 or type $2^*$ map, then at least one of the $\phi(x_i)$ is linear, and hence, at least one of the $\phi(x_i)$ has one point at infinity modulo $p$.

Let $n$ be a positive integer and assume that the statement of the theorem holds for all $p-$morphism that are a product of $n$ or less type 1 and type $2^*$ maps. Let $\phi = \phi_{n+1} \cdots \phi_2 \phi_1$, where the $\phi_i$ are type 1 and type $2^*$ maps, for $i = 1, 2, \cdots, n+1$. Let $\sigma = \phi_{n+1} \cdots \phi_3 \phi_2$ and $\sigma(x_i) = g_i$ and $\phi_1(x_i) = u_i$ for $i = 1, 2$. Then $\phi(x_i) = \sigma \phi_1(x_i) = u_i(g_1, g_2)$, for $i = 1, 2$. Each $g_i$ has one point at infinity modulo $p$ by the induction hypothesis and, if $\deg(g_1) + \deg(g_2) > 2$, each is associated with the same linear homogeneous polynomial in $k[x_1, x_2]$ by Proposition 2.16.

Case 1. If $\phi_1$ is type $2^*$, then we may assume that $u_1 = x_1$. Then $\phi(x_1) = g_1$, which has one point at infinity modulo $p$.

Case 2. If $\phi_1$ is type 1, then $u_1 = a_{11} x_1 + a_{12} x_2 + a_{13}$ and $u_2 = a_{21} x_1 + a_{22} x_2 + a_{23}$ with the $a_{ij} \in k$ and $a_{11} a_{22} - a_{12} a_{21} \neq 0$. If both $\overline{g_1}$ and $\overline{g_2}$ are linear (i.e. $\deg(\overline{g_1}) + \deg(\overline{g_2}) = 2$), then so is $\overline{u_i(g_1, g_2)}$ for each $i$ and, hence, each $\phi(x_i)$ has one point at infinity modulo $p$.

Suppose then that $\deg(\overline{g_1}) + \deg(\overline{g_2}) > 2$. Then by Proposition 2.16, for each $i = 1, 2$, $\overline{g_i}^+ = w^{e_i} H_i$, where $w_i \in k[x_1, x_2]$ is linear homogeneous, each $H_i \in k[x_1^p, x_2^p]$ is homogeneous, and the $e_i$ are positive integers not divisible by $p$. We may assume, without loss of generality, that $\deg(\overline{g_1}) \geq \deg(\overline{g_2})$. If $\deg(\overline{g_1}) > \deg(\overline{g_2})$, then $\overline{u_1(g_1, g_2)}^+ = a_{11}\overline{g_1}^+$ if $a_{11} \neq 0$ and $\overline{u_2(g_1, g_2)}^+ = a_{21}\overline{g_1}^+$ if $a_{21} \neq 0$. Since $a_{11}a_{22} - a_{12}a_{21} \neq 0$, we have $a_{11} \neq 0$ or $a_{21} \neq 0$. Hence, at least one of the $\phi(x_i)$ has one point at infinity modulo $p$.

Lastly, if $\deg(\overline{g_1}) = \deg(\overline{g_2})$, then $e_1 \equiv e_2 \pmod{p}$ and we may assume, without loss of generality, that $e_1 \geq e_2$. Then $\overline{u_1(g_1, g_2)}^+ = w^{e_2}(a_{11}w^{e_1-e_2}H_1 + a_{12}H_2)$ if $a_{11}w^{e_1-e_2}H_1 + a_{12}H_2 \neq 0$ and $\overline{u_2(g_1, g_2)}^+ = w^{e_2}(a_{21}w^{e_1-e_2}H_1 + a_{22}H_2)$ if $a_{21}w^{e_1-e_2}H_1 + a_{22}H_2 \neq 0$. Since $a_{11}a_{22} - a_{12}a_{21} \neq 0$, at least one of the two inequalities in the preceding line must hold. If $a_{11}w^{e_1-e_2}H_1 + a_{12}H_2 \neq 0$ then $a_{11}w^{e_1-e_2}H_1 + a_{12}H_2$ is homogeneous and belongs to $k[x_1^p, x_2^p]$, which implies $\phi(x_1)$ has one point at infinity modulo $p$. Likewise, if $a_{21}w^{e_1-e_2}H_1 + a_{22}H_2 \neq 0$, then $\phi(x_2)$ has one point at infinity modulo $p$. Hence, at least one of the $\phi(x_i)$ has one point at infinity modulo $p$.

**Remark 2.19.** Above we asked, if $\phi$ is a $p$-morphism of $k[x_1, x_2]$ and $\deg(\phi(x_i)) < p$ for each $i$, must $\phi(x_1)$ and $\phi(x_2)$ have one point at infinity? Since for each $i = 1, 2$, $\deg(\phi(x_i)) < p$, the number of points at infinity of $\phi(x_i)$ modulo $p$ is equal to the number of points at infinity of $\phi(x_i)$ (Definition 2.15). Hence, by Theorem 2.18, $\phi(x_1)$ and $\phi(x_2)$ each have one point at infinity.

We also asked above, if $\phi$ is a $p$-morphism of $k[x_1, x_2]$ and $\deg(\phi(x_i)) < p$ for each $i$, must $(\phi(x_1), \phi(x_2))$ be an automorphic pair? We attempt to answer this question below, but to this end it will be convenient to describe $p$-morphisms in a slightly different but equivalent way to the description given in Definition 2.2.

**Remark 2.20.** If $\phi$ is a type $2^*$ $k$-endomorphism of $k[x_1, x_2]$ defined by $x_1 \to x_1$, $x_2 \to x_2 + g(x_1, x_2^p)$ with $g(x_1, x_2) \in k[x_1, x_2]$, then $\phi = \tau \circ \rho$, where $\tau$ is the $k$-endomorphism of $k[x_1, x_2]$ defined by $x_1 \to x_1$, $x_2 \to x_2 + g(x_1, x_2^p) - g(x_1, 0)$ and $\rho$ is the type 2 $k$-automorphism of $k[x_1, x_2]$ defined by $x_1 \to x_1$, $x_2 \to x_2 + g(x_1, 0)$. Note that $x_2^p$ is a factor of $g(x_1, x_2^p) - g(x_1, 0)$, i.e. $g(x_1, x_2^p) - g(x_1, 0) = x_2^p h(x_1, x_2^p)$, for some $h(x_1, x_2) \in k[x_1, x_2]$. We will hereafter refer to such a $k$-endomorphism of $k[x_1, x_2]$ defined by $x_1 \to x_1$, $x_2 \to x_2 + x_2^p h(x_1, x_2^p)$ with $h(x_1, x_2)$ nonzero in $k[x_1, x_2]$ as a **type 3 endomorphism**. Then every type $2^*$ $k$-endomorphism $\phi$ of $k[x_1, x_2]$ is a composition $\phi = \tau \circ \rho$ where $\rho$ is a type 2 automorphism and $\tau$ is type 3. Therefore, every $p$-morphism of $k[x_1, x_2]$ is a compositional product of type 1, type 2, and type 3 maps.

**Lemma 2.21.** Let $\rho$ be a $p$-morphism. If there are two or more type 3 morphisms in the compositional product of type 1, type 2, and type 3 maps that defines $\rho$ or if one of the type 3 morphisms in the product has degree greater than or equal to $p^2$ (see Definition 1.5), then $\max\{\deg(\rho(x_1)), \deg(\rho(x_2))\} \geq p$.

**Proof.** If there are two or more type 3 morphisms in the compositional product of type 1, type 2, and type 3 maps that defines $\rho$, then by Proposition 2.9, $\deg(\rho) \geq p^2$. Thus, in either of the cases described in the lemma, we have $p^2 \leq \deg(\rho(x_1)) \cdot \deg(\rho(x_2))$ by Proposition 1.5. Hence, for some $i = 1, 2$, $\deg(\rho(x_i)) \geq p$, i.e. $\max\{\deg(\rho(x_1)), \deg(\rho(x_2))\} \geq p$.

**Lemma 2.22.** Let $f_1, f_2$ be nonzero in $k[x_1, x_2]$ and let $\widehat{f_1} = af_1 + bf_2$ and $\widehat{f_2} = cf_1 + df_2$ with $a, b, c, d \in k$ and $ad - bc \neq 0$. Then $\max\{\deg(f_1), \deg(f_2)\} = \max\{\deg(\widehat{f_1}), \deg(\widehat{f_2})\}$.

Proof. We may assume without loss of generality that $\deg(f_1) \geq \deg(f_2)$. Then $\deg(\widehat{f_i}) \leq \deg(f_1)$ for each $i = 1, 2$.

If $\deg(f_1) > \deg(f_2)$ and $a \neq 0$, then $\deg(\widehat{f_1}) = \deg(f_1)$. Otherwise, $c \neq 0$ and $\deg(\widehat{f_2}) = \deg(f_1)$. Hence, $\max\{\deg(f_1), \deg(f_2)\} = \max\{\deg(\widehat{f_1}), \deg(\widehat{f_2})\} = \deg(f_1)$.

If $\deg(f_1) = \deg(f_2)$ and $a(f_1)^+ + b(f_2)^+ \neq 0$, then $\deg(\widehat{f_1}) = \deg(f_1)$. Otherwise, $c(f_1)^+ + d(f_2)^+ \neq 0$ and $\deg(\widehat{f_2}) = \deg(f_1)$. Hence, $\max\{\deg(f_1), \deg(f_2)\} = \max\{\deg(\widehat{f_1}), \deg(\widehat{f_2})\} = \deg(f_1)$.

**Proposition 2.23.** Let $\rho$ be a $p$-morphism of $k[x_1, x_2]$. If $\rho$ is not a $k$-automorphism of $k[x_1, x_2]$, then $\max\{\deg(\rho(x_1)), \deg(\rho(x_2))\} \geq p$.

**Proof.** By Proposition 2.9 and Corollary 2.10, $\rho$ is not a $k$-automorphism of $k[x_1, x_2]$ if and only if there is at least one type 3 morphism in the compositional product of type 1, type 2, and type 3 maps that defines $\rho$. If there are two or more type 3 morphisms in the compositional product, then $\max\{\deg(\rho(x_1)), \deg(\rho(x_2))\} \geq p$ by Lemma 2.21. Thus, we are reduced to the case where there is exactly one type 3 morphism in the compositional product. Then $\rho = \sigma_1 \phi \sigma_2$ where $\phi$ is a type 3 endomorphism and where $\sigma_1$ and $\sigma_2$ are products of type 1 and type 2 maps and are therefore automorphisms.

Since $\phi$ is type 3, we may assume $\phi(x_1) = x_1$ and $\phi(x_2) = x_2 + x_2^p h(x_1, x_2^p)$ with $h(x_1, x_2)$ a nonzero element of $k[x_1, x_2]$. If $\deg_{x_2}(h(x_1, x_2)) \geq 1$, then $\deg_p(\phi) \geq 2$ (see Definition 2.7), which implies $\deg(\phi) \geq p^2$ by Proposition 2.9, which by Lemma 2.21 yields $\max\{\deg(\rho(x_1)), \deg(\rho(x_2))\} \geq p$. Hence, to complete the proof, we may assume that $\deg_{x_2}(h(x_1, x_2)) = 0$, i.e. $\phi(x_2) = x_2 + x_2^p h(x_1)$ with $h(x_1)$ a nonzero element of $k[x_1]$.

We will proceed by induction on the number of type 1 and type 2 non-identity automorphisms in the compositional product that defines $\sigma_2$. If there are no such automorphisms (i.e. $\sigma_2$ is the identity map) then $\rho = \sigma_1 \phi$. For $i = 1, 2$, let $u_i = \sigma_1(x_i)$. Then $\rho(x_1) = u_1$ and $\rho(x_2) = u_2 + u_2^p h(u_1)$. Since $\deg(u_2(x_1, x_2)) > 0$, we have $\deg(\rho(x_2)) > p$. Hence, $\max\{\deg(\rho(x_1)), \deg(\rho(x_2))\} \geq p$.

Assume now that if $\sigma_2$ is a product of $n$ or less non-identity type 1 and type 2 maps then $\max\{\deg(\rho(x_1)), \deg(\rho(x_2))\} \geq p$. Let $\sigma_2$ be such a product and $\tau$ be a type 1 or type 2 non-identity automorphism of $k[x_1, x_2]$. We will complete the proof by showing that $\max\{\deg(\rho\tau(x_1)), \deg(\rho\tau(x_2))\} \geq p$. Without loss of generality, we may assume $\deg(\rho(x_1)) \geq \deg(\rho(x_2))$. Then $\deg(\rho(x_1)) \geq p$ by the induction hypothesis. For all but at most one $a \in k$, $\deg(\rho(x_2) + a\rho(x_1)) = \deg(\rho(x_1))$; i.e. $\deg(\rho(x_2 + ax_1)) = \deg(\rho(x_1))$. Choose such an $a$. Then by Lemma 2.22, $\max\{\deg(\rho(x_1)), \deg(\rho(x_2 + ax_1))\} = \max\{\deg(\rho(x_1)), \deg(\rho(x_2))\}$. Hence, if we replace $x_2$ by $x_2 + ax_1$, we may assume $\deg(\rho(x_1)) = \deg(\rho(x_2)) \geq p$. If $\tau$ is type 1, then $\tau(x_1) = a_{11}x_1 + a_{12}x_2 + a_{13}$ and $\tau(x_2) = a_{21}x_1 + a_{22}x_2 + a_{23}$ with the $a_{ij} \in k$ and $a_{11}a_{22} - a_{12}a_{21} \neq 0$. Then $\rho\tau(x_1) = a_{11}\rho(x_1) + a_{12}\rho(x_2) + a_{13}$ and $\rho\tau(x_2) = a_{21}\rho(x_1) + a_{22}\rho(x_2) + a_{23}$. Then by Lemma 2.22, $\max\{\deg(\rho\tau(x_1)), \deg(\rho\tau(x_2))\} = \max\{\deg(\rho(x_1)), \deg(\rho(x_2))\} \geq p$. If $\tau$ is type 2, then we may assume that $\tau(x_1) = x_1$ and $\tau(x_2) = x_2 + g(x_1)$ where $g \in k[x_1]$. Then $\rho\tau(x_1) = \rho(x_1)$ and $\rho\tau(x_2) = \rho(x_2) + g(\rho(x_1))$. Hence, $\max\{\deg(\rho\tau(x_1)), \deg(\rho\tau(x_2))\} \geq \deg(\rho(x_1)) \geq p$.

**Section 3. Other Methods of Generating Jacobian Pairs**

In this section and the next we investigate the existence of other possible methods of generating Jacobian pairs with a view to testing them regarding the Separable Jacobian Conjecture and the Low Degree Jacobian Conjectures mentioned above. As we have seen, if there are no others, then the Separable Jacobian Conjecture in characteristic $p > 0$ is true. Of course, when we obtain seemingly alternative methods of generating Jacobian pairs, the first question we should ask is, are there pairs generated by one method that the others do not? When I occasionally investigated characteristic $p$ approaches to the Jacobian Conjecture years ago, long before I considered $p$-morphisms, I would simply create examples by taking an existing Jacobian pair $(f, g)$ and then replace $g$ by $\tilde{g} = g + \lambda_0 + \lambda_1 f + \cdots + \lambda_r f^r$ where $r$ a positive integer and the $\lambda_i$ are in $k[g^p]$, to obtain a new Jacobian pair $(f, \tilde{g})$. And then repeat the process with $(f, \tilde{g})$, etcetera. We formalize this approach below and seek to determine if this method and the $p$-morphism approach produce different classes of Jacobian pairs.

**Definition 3.1.** Let $A = k[x_1, x_2]$. If $\phi$ is a $k$-endomorphism of $A$ defined by $\phi(x_i) = f_i(x_1, x_2)$ where $f_i \in A$ for $i = 1, 2$, we let $\bar{\phi}: A \times A \to A \times A$ denote the map defined by $\bar{\phi}(g_1, g_2) = (f_1(g_1, g_2), f_2(g_1, g_2))$, for all $(g_1, g_2) \in A \times A$. We refer to $\bar{\phi}$ as the **induced map on $A \times A$ by $\phi$**. Note that $\bar{\phi}(x_1, x_2) = (\phi(x_1), \phi(x_2))$ but in general, $\bar{\phi}(g_1, g_2) \neq (\phi(g_1), \phi(g_2))$ for $(g_1, g_2) \in A \times A$.

**Remark 3.2.** If $\phi$ is a $k$-endomorphism of $A = k[x_1, x_2]$ defined by $\phi(x_i) = f_i(x_1, x_2)$ where $f_i \in A$ for $i = 1, 2$ and if $(f_1, f_2)$ is a Jacobian pair, then for all Jacobian pairs $(g_1, g_2) \in A \times A$, $\bar{\phi}(g_1, g_2)$ will be a Jacobian pair by the Multivariable Chain Rule. Thus, if $\phi$ is a $p$-morphism, then $\bar{\phi}(g_1, g_2)$ will be a Jacobian pair for all Jacobian pairs $(g_1, g_2) \in A \times A$.

**Definition 3.3** (The **Induced Map Method** of Generating Jacobian Pairs). For each $n = 1, 2, 3, \cdots$, let $\phi_n$ be a type 1 or $2^*$ morphism on $A = k[x_1, x_2]$ and let $\bar{\phi}_n$ be the corresponding induced map on $A \times A$. For each positive integer $n$, let $(f_{n1}, f_{n2}) = \bar{\phi}_n \cdots \bar{\phi}_3 \bar{\phi}_2 \bar{\phi}_1 (x_1, x_2)$. Then $(f_{n1}, f_{n2})$ will be a Jacobian pair for each $n$ by Remark 3.2.

**Example 3.4.** As an example of the Induced Map Method, let $\phi_1$ be the type $2^*$ morphism on $k[x_1, x_2]$ defined by $\phi_1(x_1) = x_1$, $\phi_1(x_2) = x_2 + x_1^3$, let $\phi_2$ be the type $2^*$ morphism on $k[x_1, x_2]$ defined by $\phi_2(x_1) = x_1 + x_2^5 x_1^p$, $\phi_2(x_2) = x_2$, and let $\phi_3$ be the type $2^*$ morphism on $k[x_1, x_2]$ defined by $\phi_3(x_1) = x_1$, $\phi_3(x_2) = x_2 - x_1^7$. Then by the Induced Map Method we obtain, in order, the Jacobian pairs, $(f_{11}, f_{12}) = (x_1, x_2 + x_1^3)$, $(f_{21}, f_{22}) = (x_1 +$

$(x_2 + x_1{}^3)^5 x_1{}^p, x_2 + x_1{}^3)$ and $(f_{31}, f_{32}) = (x_1 + (x_2 + x_1{}^3)^5 x_1{}^p, x_2 + x_1{}^3 - (x_1 + (x_2 + x_1{}^3)^5 x_1{}^p)^7)$.

**Remark 3.5.** I have found the Induced Map Method an efficient way to generate Jacobian pairs. Eventually, I realized that this method of generating Jacobian pairs and that by using $p-$morphisms are equivalent, in the sense that any example produced by one can also be produced by the other, which is a corollary of the next result.

**Proposition 3.6.** Let $\sigma$ and $\tau$ be $p-$morphisms on $A = k[x_1, x_2]$. Let $\bar{\sigma}$ and $\bar{\tau}$ be the corresponding induced maps on $A \times A$. Then $(\bar{\tau} \circ \bar{\sigma})(x_1, x_2) = \big((\sigma \circ \tau)(x_1), (\sigma \circ \tau)(x_2)\big)$.

**Proof.** We have $\sigma(x_i) = f_i(x_1, x_2)$ and $\tau(x_i) = g_i(x_1, x_2)$ for some $f_i, g_i \in k[x_1, x_2]$, for $i = 1, 2$. Then $(\bar{\tau} \circ \bar{\sigma})(x_1, x_2) = \bar{\tau}\big(f_1(x_1, x), f_2(x_1, x_2)\big) = \big(g_1(f_1(x_1, x_2)), g_2(f_2(x_1, x_2))\big)$, while $(\sigma \circ \tau)(x_i) = \sigma(\tau(x_i)) = \sigma(g_i(x_1, x_2)) = g_i(f_1(x_1, x_2), f_2(x_1, x_2))$ for $i = 1, 2$.

**Corollary 3.7.** Let $\sigma_1, \sigma_2, \cdots, \sigma_n$ be $p-$morphisms on $A = k[x_1, x_2]$ and $\bar{\sigma}_1, \bar{\sigma}_2, \cdots, \bar{\sigma}_n$ be the corresponding induced maps on $A \times A$. Then,

$$(\bar{\sigma}_n \cdots \bar{\sigma}_2 \bar{\sigma}_1)(x_1, x_2) = \big((\sigma_1 \sigma_2, \cdots \sigma_n)(x_1), ((\sigma_1 \sigma_2, \cdots \sigma_n))(x_2)\big).$$

**Proof** by induction on $n$. The $n = 1$ case is covered in Definition 2.1. Assume the result for the $n - 1$ case. We have $\sigma_n(x_i) = f_i$ for some $f_i \in k[x_1, x_2]$, $i = 1, 2$. Then by the induction hypothesis,

$$(\bar{\sigma}_n \cdots \bar{\sigma}_2 \bar{\sigma}_1)(x_1, x_2) = \bar{\sigma}_n\big((\bar{\sigma}_{n-1} \cdots \bar{\sigma}_2 \bar{\sigma}_1)(x_1, x_2)\big) =$$
$$\bar{\sigma}_n\big((\sigma_1 \sigma_2, \cdots \sigma_{n-1})(x_1), ((\sigma_1 \sigma_2, \cdots \sigma_{n-1}))(x_2)\big) =$$
$$\big((\sigma_1 \sigma_2, \cdots \sigma_{n-1})(f_1(x_1, x_2)), ((\sigma_1 \sigma_2, \cdots \sigma_{n-1}))(f_2(x_1, x_2))\big) =$$
$$\big((\sigma_1 \sigma_2, \cdots \sigma_n)(x_1), ((\sigma_1 \sigma_2, \cdots \sigma_n))(x_2)\big).$$

**Remark 3.8.** From Corollary 3.7 we see that the induced map method of generating Jacobian pairs and the $p-$morphism approach produce the same sets of Jacobian pairs. Even so, for producing actual examples of Jacobian pairs I have found the induced map method easier to implement.

**Section 4. Jacobian Pairs Generated by Functions of One Variable.**

Jacobian pairs that are generated by $p-$morphisms have one point at infinity modulo $p$ by Proposition 2.18. Thus, in search of Jacobian pairs that are not produced by $p-$morphisms, it makes sense to try to construct Jacobian pairs with two or more points at infinity modulo $p$ (see Theorem 2.18). In Example 4.4 below, a special case of which we first introduced in [12](page 51), we obtain such Jacobian pairs.

**Example 4.1.** Let $f_1 = \alpha x_1^a x_2^b + x_1$, $f_2 = \beta x_1^c x_2^d + x_2$ with $\alpha, \beta \in k^*$, and $a, b, c, d$ positive integers such that $abcd \not\equiv 0 \pmod{p}$. Then,

$$\frac{\partial(f_1, f_2)}{\partial(x_1, x_2)} = (ad - bc)\alpha\beta x_1^{a+c-1} x_2^{b+d-1} + a\alpha x_1^{a-1} x_2^b + d\beta x^c x_2^{d-1} + 1.$$

Then, $\partial(f_1, f_2)/\partial(x_1, x_2) = 1$ if and only if $ad \equiv bc, a - 1 = c, b = d - 1, a\alpha + d\beta = 0$. A pair $f_1, f_2$ as above will then have the form,

$$f_1 = x_1(\alpha x_1^{a-1} x_2^{mp-a} + 1) \text{ and } f_2 = x_2\left(\frac{a\alpha}{a-1} x_1^{a-1} x_2^{mp-a} + 1\right),$$

where $a \not\equiv 1$ and $m, n, r$ are positive integers with $mp > a$. A straightforward computation will confirm that $(f_1, f_2)$ as above will be a Jacobian pair.

A Jacobian pair of this kind with $\deg(f_1) + \deg(f_2)$ a minimum will have the form,

$$f_1 = x_1(\alpha x_1^{a-1} x_2^{p-a} + 1) \text{ and } f_2 = x_2\left(\frac{a}{a-1} \alpha x_1^{a-1} x_2^{p-a} + 1\right),$$

where $1 < a < p$.

Note that the Jacobian pairs in Example 4.1 have the form, $f_1 = x_1 h_1(u)$ where $h_1, h_2 \in k[x]$ are linear polynomials and $u = x_1^{a-1} x_2^{mp-a}$. In the rest of this section we seek to find to what extent this example generalizes.

**Proposition 4.2.** *Let $s \leq t \leq d$ be nonnegative integers and let $\alpha_j \in k$ for $s \leq j \leq t$. Let $u = \sum_{j=s}^{t} \alpha_j x_1^j x_2^{d-j} \in k[x_1, x_2]$ with $\alpha_s \alpha_t \neq 0$. For each $i = 1, 2$, let $f_i = x_i h_i(u)$ where $h_i \in k[x]$. If $\partial(f_1, f_2)/\partial(x_1, x_2) \in k^*$, then $[k(x_1, x_2): k(f_1, f_2)] = t(deg(h_1)) + (d - s)(deg(h_2)) + 1$.*

**Proof.** By Remark 2.3 and since $f_i/h_i(u) = x_i$, for $i = 1, 2$, we have $k(x_1, x_2)$ is a primitive algebraic extension of $k(f_1, f_2)$ generated by $u$. Hence, $[k(x_1, x_2): k(f_1, f_2)]$ equals the degree of $u$ over $k(f_1, f_2)$.

We have,

$$\sum_{j=s}^{t} \alpha_j (h_1(u))^{d-j} (h_2(u))^j (f_1)^j (f_2)^{d-j} = \sum_{j=s}^{t} \alpha_j x_1^j x_2^{d-j} (h_1(u))^d (h_2(u))^d$$

$$= u(h_1(u) h_2(u))^d.$$

Hence, $u$ is a root of the polynomial $H(T)$ in $k[f_1, f_2][T]$ given by,

$$H(T) = T(h_1(T) h_2(T))^d - \sum_{j=s}^{t} \alpha_j (h_1(T))^{d-j} (h_2(T))^j (f_1)^j (f_2)^{d-j}.$$

Dividing both sides of this equality by $(h_1(T))^{d-t}(h_2(T))^s$, we obtain $u$ is a root of the polynomial $M(T)$ in $k[f_1, f_2][T]$ given by,

$$M(T) = T(h_1(T))^t(h_2(T))^{d-s} - \sum_{j=s}^{t} \alpha_j (h_1(T))^{t-j}(h_2(T))^{j-s}(f_1)^j(f_2)^{d-j}.$$

Since $\deg_T(M(T)) = t(\deg(h_1)) + (d-s)(\deg(h_2)) + 1$, we will be done if we show that $M(T)$ is irreducible in $k(f_1, f_2)[T]$. To this end let $X_1 = f_1$ and $X_2 = \frac{f_2}{f_1}$. Then, as a polynomial in $k[X_2, T][X_1]$, we have,

$$M(T) = T(h_1(T))^t(h_2(T))^{d-s} - \sum_{j=s}^{t} \alpha_j (h_1(T))^{t-j}(h_2(T))^{j-s}(X_2)^{d-j}(X_1)^d.$$

Since $\partial(f_1, f_2)/\partial(x_1, x_2) \in k^*$, we have $h_i(0) \neq 0$ for each $i$. By Remark 2.3, $h_1(T)$ and $h_2(T)$ have no common factors in $k[T]$. It follows that $T$, $h_1(T)$ and $h_2(T)$ are pairwise relatively prime in $k[T]$. Thus, the greatest common divisor of $T(h_1(T))^t(h_2(T))^{d-s}$ and $\sum_{j=s}^t \alpha_j (h_1(T))^{t-j}(h_2(T))^{j-s}(X_2)^{d-j}$ in $k[X_2, T]$ is 1. It follows that the content of $M$ as a polynomial in $X_1$ with coefficients in $k[X_2, T]$ equals 1. Since $T$ divides $T(h_1(T))^t(h_2(T))^{d-s}$ in $k[X_2, T]$ but $T^2$ does not and since $T$ does not divide $\sum_{j=s}^t \alpha_j (h_1(T))^{t-j}(h_2(T))^{j-s}(X_2)^{d-j}$ in $k[X_2, T]$, we have that $M$ is irreducible in $k[X_1, X_2, T]$ by Eisenstein's Lemma. Then by Gauss's Lemma, $M$ is irreducible as a polynomial in $k(X_1, X_2)[T] = k(f_1, f_2)[T]$, which proves the result.

**Corollary 4.3.** *Let* $u = x_1^a x_2^b$ *with* $a, b$ *nonnegative integers and for* $i = 1, 2$, *let* $f_i = x_i h_i(u)$ *where* $h_i \in k[x]$. *If* $\partial(f_1, f_2)/\partial(x_1, x_2) \in k^*$, *then* $[k(x_1, x_2): k(f_1, f_2)] = a \cdot deg(h_1) + b \cdot deg(h_2) + 1$.

**Proof.** In the notation of Proposition 4.2 we have $s = t = a$ and $d = a + b$, from which the result follows.

**Example 4.4.** Let $u = x_1^a x_2^b$ with $a, b$ positive integers and for $i = 1, 2$, let $f_i = x_i h_i(u)$ where $h_i \in k[x]$. We seek to find Jacobian pairs of the form $(f_1, f_2)$ with $ab \not\equiv 0 \pmod{p}$, for otherwise, one of the $f_i$ will have one point at infinity modulo $p$ (so will the other by Lemma 2.17) and we want the $f_i$ to each have two. We will also assume $\deg(h_i) > 0$ for each $i$ for the same reason. We then have,

$$\frac{\partial(f_1, f_2)}{\partial(x_1, x_2)} = \det \begin{bmatrix} h_1(u) + x_1 h_1'(u) u_{x_1} & x_1 h_1'(u) u_{x_2} \\ x_2 h_2'(u) u_{x_1} & h_2(u) + x_2 h_2'(u) u_{x_2} \end{bmatrix}$$

$$= h_1(u)h_2(u) + h_1'(u)h_2(u)x_1 u_{x_1} + h_1(u)h_2'(u)x_2 u_{x_2}$$

$$= h_1(u)h_2(u) + au \cdot h_1'(u)h_2(u) + bu \cdot h_1(u)h_2'(u).$$

It follows that,

(4.4.1)

$$\partial(f_1, f_2)/\partial(x_1, x_2) \in k^* \text{ if and only if } h_1(u)h_2(u) + auh_1'(u)h_2(u) + buh_1(u)h_2'(u) \in k^*$$

**Case 1.** $a = b = 1$.

Then $\partial(f_1, f_2)/\partial(x_1, x_2) \in k^*$ if and only if $h_1(u)h_2(u) + u \cdot h_1'(u)h_2(u) + u \cdot h_1(u)h_2'(u) \in k^*$, which is equivalent to $\frac{d}{du}(uh_1(u)h_2(u)) \in k^*$; i.e. $uh_1(u)h_2(u) = \alpha u + g(u^p)$, for some $\alpha \in k^*$ and some $g \in k[x]$ with $g(0) = 0$.

Hence, in this case, $(f_1, f_2)$ will be a Jacobian pair if and only if $h_1(u)h_2(u) = \alpha + \frac{g(u^p)}{u}$, for some $\alpha \in k^*$ and $g \in k[x]$ with $g(0) = 0$. By Proposition 21, $[k(x_1, x_2): k(f_1, f_2)] = \deg(h_1) + \deg(h_2) + 1 = p \cdot \deg(g) \equiv 0 \pmod{p}$.

Also, $\deg(f_1 f_2) = 2p \cdot \deg(g) \geq 2p$, which implies at least one of the $f_i$ has degree greater than or equal to $p$.

**Case 2.** $a + b > 2$.

We have $h_1(u) = \sum_{i=0}^m \alpha_i u^i$ and $h_2(u) = \sum_{j=0}^n \beta_j u^j$ for some $\alpha_i, \beta_j \in k$ and positive integers $m$ and $n$ with $\alpha_m \beta_n \neq 0$. Without loss of generality, we may assume $m \geq n$. Then by (4.4.1), $\partial(f_1, f_2)/\partial(x_1, x_2) \in k^*$ if and only if

$$\sum_{i=0}^m \sum_{j=0}^n (1 + bj + ai) \alpha_i \beta_j u^{i+j} \in k^*.$$

In particular, the coefficient of $u^{m+n}$ must be 0, i.e. $1 + am + bn \equiv 0$. By Corollary 4.3, we obtain $[k(x_1, x_2): k(f_1, f_2)] = a \cdot \deg(h_1) + b \cdot \deg(h_2) + 1 = am + bn + 1 \equiv 0$.

We also obtain $\deg(f_1) + \deg(f_2) = (a + b)(m + n) + 2 = (am + bn + 1) + (an + bm + 1) \geq (am + bn + 1) + (an + bn + 1) = (am + bn + 1) + \deg(f_2)$. Thus, $\deg(f_1) \geq am + bn + 1 \equiv 0$. It follows that $\deg(f_1) = ps$, for some positive integer $s$. Hence, $\deg(f_1) \geq ps \geq p$. Therefore, for Jacobian pairs of the form described in this case, we have $[k(x_1, x_2): k(f_1, f_2)] \equiv 0 \pmod{p}$ and at least one of the $f_i$ has degree greater than or equal to $p$.

To produce a concrete example of a Jacobian pair as described in case 2, let $m = 2$ and $n = 1$. Then $h_1(u) = \alpha_0 + \alpha_1 u + \alpha_2 u^2$ and $h_2(u) = \beta_0 + \beta_1 u$. Then $\frac{\partial(f_1, f_2)}{\partial(x_1, x_2)} \in k^*$ if and only if $\alpha_0 \beta_0 \in k^*$, $\alpha_1 \beta_0 (1 + a) + \alpha_0 \beta_1 (1 + b) = 0$, $\alpha_2 \beta_0 (1 + 2a) + \alpha_1 \beta_1 (1 + a + b) = 0$, and $\alpha_2 \beta_1 (1 + 2a + b) = 0$. We may assume that $\alpha_0 = \beta_0 = 1$. Solving this system, we obtain the following conditions on the $h_i$ for $(f_1, f_2)$ to be a Jacobian pair: $b \equiv -1 - 2a$, $\beta_1 = \left(\frac{1+a}{2a}\right)\alpha_1$, $\alpha_2 = \frac{a+1}{2a(2a+1)}\alpha_1^2$, $a \not\equiv 0$, $a \not\equiv -1$, $a \not\equiv -\frac{1}{2}$, $\alpha_1 \neq 0$, $p > 2$. It follows that if $p > 2$ and $s$ is a positive integer such that $sp > -1 - 2a$, then $f_1 = x_1\left(1 + \alpha_1 u + \frac{a+1}{2a(2a+1)}\alpha_1^2 u^2\right)$ and $f_2 = x_2\left(1 + \left(\frac{1+a}{2a}\right)\alpha_1 u\right)$ with $u = x_1{}^a x_2{}^{sp-1-2a}$ is a Jacobian pair with $\deg(f_1) = 2(sp - a) - 1$ and $\deg(f_2) = sp - a$.

Note that when $s = 1$, we get $u = x_1{}^a x_2{}^{p-1-2a}$. Hence, $a < \frac{p-1}{2}$. Then $\deg(f_1) = 2(p - a) - 1 \geq p$ and $\deg(f_2) = p - a < p$. Therefore, $f_1$ and $f_2$ have two points at infinity, but only $f_2$ has degree less than $p$.

The following proposition summarizes the main conclusions of Example 4.4.

**Proposition 4.5.** *Let $u = x_1{}^a x_2{}^b$ with $a, b$ positive integers such that $ab \not\equiv 0 \pmod{p}$ and for $i = 1, 2$, let $f_i = x_i h_i(u)$ where $h_i \in k[x]$. Assume $\deg(h_i) > 0$ for each $i$. Then $\partial(f_1, f_2)/\partial(x_1, x_2) \in k^*$ if and only if $h_1(u)h_2(u) + auh_1'(u)h_2(u) + buh_1(u)h_2'(u) \in k^*$. If such is the case, then $[k(x_1, x_2): k(f_1, f_2)] \equiv 0 \pmod{p}$ and at least one of the $f_i$ has degree greater than or equal to $p$.*

At this point, it is natural to investigate if there exists Jacobian pairs of the form described in Proposition 4.5 but with $u \in k[x_1, x_2]$ a homogeneous form with more than two pairwise relatively prime linear factors. This is the subject of the next example, but first we'll need a preliminary lemma.

**Lemma 4.6.** *Let $d$ be a positive integer and $a_i \in k$ for $0 \leq i \leq d$. Let $u = \sum_{i=0}^{d} a_i x_1{}^i x_2{}^{d-i} \in k[x_1, x_2]$ with $u \neq 0$. Let $f_1 = x_1 h_1(u), f_2 = x_2 h_2(u)$, where $h_1, h_2 \in k[x]$. Let the degrees of $h_1$ and $h_2$ be $m$ and $n$, respectively. If $\deg(\partial(f_1, f_2)/\partial(x_1, x_2)) < d(m + n)$, then $((m - n)i + 1 + nd)a_i = 0$ for each $i$ such that $0 \leq i \leq d$.*

**Proof.** We have $h_1 = \sum_{j=0}^{m} \alpha_j t^j$ and $h_2 = \sum_{j=0}^{n} \beta_j t^j$ for some $\alpha_j, \beta_j \in k$ with $\alpha_m \beta_n \neq 0$. Then $\partial(f_1, f_2)/\partial(x_1, x_2) = h_1(u)h_2(u) + h_1'(u)h_2(u)x_1 u_{x_1} + h_1(u)h_2'(u)x_2 u_{x_2}$. We also have $\deg(\partial(f_1, f_2)/\partial(x_1, x_2)) \leq \deg(f_1) + \deg(f_2) - 2 = d(m + n)$. Hence, $(\partial(f_1, f_2)/\partial(x_1, x_2)) = \sum_{i=0}^{m+n} H_i$, where each $H_i$ belongs to the $k$−space generated by the monomials in $k[x_1, x_2]$ of degree $i$. By a straightforward computation we have,

$$H_{m+n} = \alpha_m \beta_n u^{m+n-1}\left(u + m\sum_{i=0}^{d} ia_i x_1{}^i x_2{}^{d-i} + n\sum_{i=0}^{d}(d-i)a_i x_1{}^i x_2{}^{d-i}\right)$$

$$= \alpha_m \beta_n u^{m+n-1}\sum_{i=0}^{d}(1 + mi + n(d - i))a_i x_1{}^i x_2{}^{d-i}.$$

It follows that if $\deg(\partial(f_1, f_2)/\partial(x_1, x_2)) < d(m + n)$, then,

$$\sum_{i=0}^{d}(1 + mi + n(d - i))a_i x_1{}^i x_2{}^{d-i} = 0.$$

Hence, if $\deg(\partial(f_1, f_2)/\partial(x_1, x_2)) < d(m+n)$, then $((m-n)i + 1 + nd)a_i = 0$, for each $i$ such that $0 \leq i \leq d$.

**Example 4.7.** Let $d$ be a positive integer and $a_i \in k$ for $0 \leq i \leq d$. Let $u = \sum_{i=0}^{d} a_i x_1^i x_2^{d-i} \in k[x_1, x_2]$ be nonzero. Let $f_1 = x_1 h_1(u), f_2 = x_2 h_2(u)$, where $h_1, h_2 \in k[x]$. Let the degrees of $h_1$ and $h_2$ be positive integers $m$ and $n$, respectively. Assume $\partial(f_1, f_2)/\partial(x_1, x_2) \in k^*$.

**Case 1.** Assume $m \not\equiv n \pmod{p}$. Then by Lemma 4.6 we have for each $i = 0, 1, \cdots, d$, $a_i = 0$ or $i \equiv \frac{1+nd}{n-m}$. It follows that $u = x_1^a x_2^b (h(x_1, x_2))^p$, for some nonnegative integers $a, b$, with $a \equiv \frac{1+nd}{n-m}$, $h(x_1, x_2) \in k[x_1, x_2]$ is homogeneous, and $a + b + p\deg(h) = d$. Then by the same reasoning used in Example 4.4, we again have $h_1(u)h_2(u) + auh_1'(u)h_2(u) + buh_1(u)h_2'(u) \in k^*$. Also note that $a \equiv \frac{1+nd}{n-m}$ if and only if $an - am \equiv 1 + nd$, i.e. $n(d-a) + ma + 1 \equiv 0$; i.e. $ma + nb + 1 \equiv 0 \pmod{p}$, as is the case for the Jacobian pairs in Example 4.4. Therefore, the Jacobian pairs of Example 4.4 are a special case of the Jacobian pairs of the form considered here.

In the case under consideration, the Jacobian pairs satisfy $s \equiv t \equiv a$ and $d \equiv a + b$ in the notation of the statement of Proposition 4.2. Hence, $[k(x_1, x_2): k(f_1, f_2)] \equiv 0 \pmod{p}$. Further, if $\deg(h) > 0$, then both of the $f_i$ have degree greater than or equal to $p$, If $\deg(h) = 0$, then at least one of the $f_i$ has degree greater than or equal to $p$ by Proposition 4.5. Finally, note that $ab \not\equiv 0 \pmod{p}$. For if $a \equiv 0 \pmod{p}$, then $u \in k[x_1^p, x_2]$ and is homogeneous of positive degree, and since $\deg(h_2) = n > 0$, $f_2 \in k[x_1^p, x_2]$ and has degree greater than 2, which precludes $\partial(f_1, f_2)/\partial(x_1, x_2) \in k^*$.

**Case 2.** Assume $m \equiv n \pmod{p}$. Then by Lemma 4.6, $1 + nd \equiv 0 \pmod{p}$. Then $x_1 u_{x_1} = \sum_{i=0}^{d} i a_i x_1^i x_2^{d-i}$ and $x_2 u_{x_2} = \sum_{i=0}^{d} (d-i) a_i x_1^i x_2^{d-i}$ and they are linearly dependent over $k$ if and only if $(i(d-j) - j(d-i))a_i a_j = 0$ for all $i, j = 1, 2, \cdots, d$; i.e. $d(i-j)a_i a_j = 0$ for all $i, j$; i.e. $(i-j)a_i a_j = 0$ for all $i, j$ since $d \not\equiv 0 \pmod{p}$. Thus, $x_1 u_{x_1}$ and $x_2 u_{x_2}$ are linearly dependent over $k$ if and only if $u = x^a y^b (h(x, y))^p$ for some nonnegative integers $a, b$ and homogeneous $h(x_1, x_2) \in k[x_1, x_2]$. Since this case is covered in Case 1, we will assume from here on that $x_1 u_{x_1}$ and $x_2 u_{x_2}$ are linearly independent over $k$.

We have $h_1 = \sum_{j=0}^{m} \alpha_j t^j$ and $h_2 = \sum_{j=0}^{n} \beta_j t^j$ for some nonnegative integers $m$ and $n$ and $\alpha_j, \beta_j \in k$ with $\alpha_m \beta_n \neq 0$. We also have $\alpha_0 \beta_0 \neq 0$ since $\partial(f_1, f_2)/\partial(x_1, x_2) \in k^*$. For all integers $i > m, j > n$, and $i < 0, j < 0$, define $\alpha_i = \beta_j = 0$. We have,

$$\partial(f_1, f_2)/\partial(x_1, x_2) = h_1(u)h_2(u) + h_1'(u)h_2(u)x_1 u_{x_1} + h_1(u)h_2'(u)x_2 u_{x_2}$$

$$= \sum_{i=0}^{m} \sum_{j=0}^{n} \left( \alpha_i \beta_j u^{i+j} + x_1 u_{x_1} i \alpha_i \beta_j u^{i+j-1} + x_j u_{x_j} j \alpha_i \beta_j u^{i+j-1} \right)$$

$$= \sum_{i=0}^{m} \sum_{j=0}^{n} \alpha_i \beta_j u^{i+j-1} (u + i x_1 u_{x_1} + j x_2 u_{x_2})$$

$$= \sum_{i=0}^{m} \sum_{j=0}^{n} \alpha_i \beta_j u^{i+j-1} \left( \frac{x_1 u_{x_1} + x_2 u_{x_2}}{d} + i x_1 u_{x_1} + j x_2 u_{x_2} \right)$$

$$= \sum_{i=0}^{m} \sum_{j=0}^{n} \alpha_i \beta_j u^{i+j-1} \left( -m(x_1 u_{x_1} + x_2 u_{x_2}) + i x_1 u_{x_1} + j x_2 u_{x_2} \right)$$

$$= \sum_{i=0}^{m} \sum_{j=0}^{n} \alpha_i \beta_j u^{i+j-1} \left( (i-m) x_1 u_{x_1} + (j-m) x_2 u_{x_2} \right)$$

$$= \sum_{r=0}^{m+n} \sum_{i=0}^{r} \alpha_i \beta_{r-i} u^{r-1} \left( (i-m) x_1 u_{x_1} + (r - (i+m)) x_2 u_{x_2} \right)$$

$$= \sum_{r=0}^{m+n} \sum_{i=0}^{r} \alpha_i \beta_{r-i} \left( (i-m) x_1 u_{x_1} + (r - (i+m)) x_2 u_{x_2} \right) u^{r-1}$$

Hence, $\partial(f_1, f_2)/\partial(x_1, x_2) \in k^*$ if and only if for each $r = 1, 2, \cdots, m+n$,

$$\sum_{i=0}^{r} \alpha_i \beta_{r-i} \left( (i-m) x_1 u_{x_1} + (r - (i+m)) x_2 u_{x_2} \right) = 0.$$

Since $x_1 u_{x_1}$ and $x_2 u_{x_2}$ are linearly independent over $k$, $\partial(f_1, f_2)/\partial(x_1, x_2) \in k^*$ if and only if for each $r = 1, 2, \cdots, m+n$,

$$(4.7.1a) \sum_{i=0}^{r} (i-m) \alpha_i \beta_{r-i} = 0 \text{ and } (4.7.1b) \sum_{i=0}^{r} (r - (i+m)) \alpha_i \beta_{r-i} = 0.$$

**Claim.** If $\partial(f_1, f_2)/\partial(x_1, x_2) \in k^*$, then $\alpha_{m-i} = \beta_{n-i} = 0$ for each integer $i \not\equiv 0 \pmod{p}$.

**Proof of Claim.** Since by the above definition $\alpha_i = \beta_i = 0$ for all integers $i < 0$, it is enough to prove the claim for all nonnegative integers $i$. We will proceed by induction on $i$. Thus, our initial case is when $i = 1$. Then by definition, in equation (4.7.1a) with $r = m + n - 1$ we will have $\alpha_i = 0$ for $i \geq m+1$ and $\beta_{r-i} = 0$ for $i \leq m-2$ since $i \leq m-2$ implies $r - i \geq n+1$. Then from (4.7.1.a) we have $-\alpha_{m-1} \beta_n = 0$ and from (4.7.1.b) we have $-\alpha_m \beta_{n-1} = 0$. Since $\alpha_m \beta_n \neq 0$, $\alpha_{m-1} = \beta_{n-1} = 0$.

To complete the induction, assume $s$ is a positive integer with $s \not\equiv 0 \pmod{p}$ and assume $\alpha_{m-i} = \beta_{n-i} = 0$ for each $i \not\equiv 0 \pmod{p}$ with $i < s$. We will show that $\alpha_{m-s} = \beta_{n-s} = 0$. By the above definitions we are done if $s \geq m + n$. Thus, we may assume that $s < m + n$. By (4.7.1a) with $r = m + n - s$, we have

$$(4.7.2) \sum_{i=0}^{m-(s+1)} (i-m) \alpha_i \beta_{r-i} + \sum_{i=m-s}^{m} (i-m) \alpha_i \beta_{r-i} + \sum_{i=m+1}^{m+n-s} (i-m) \alpha_i \beta_{r-i} = 0.$$

If $0 \leq i \leq m - (s+1)$, then $r - i = m + n - s - i > n$, which implies $\beta_{r-i} = 0$ by the above definition, which implies $\sum_{i=0}^{m-(s+1)} (i-m) \alpha_i \beta_{r-i} = 0$. If $m + 1 \leq i \leq m + n - s$, then $\alpha_i = 0$ by the above definition, which implies $\sum_{i=m+1}^{m+n-s}(i-m)\alpha_i\beta_{r-i} = 0$. Hence, (4.7.1a) reduces to the equation,

$$(4.7.2a) \sum_{i=m-s}^{m-1}(i-m)\alpha_i\beta_{r-i} = 0.$$

Likewise, (4.7.1b) reduces to the equation,

$$(4.7.2b) \sum_{i=m-s}^{m}(r-(i+m))\alpha_i\beta_{r-i} = 0$$

Letting $j = i - (m-s)$, these two equations become,

$$(4.7.3a) \sum_{j=0}^{s-1}(j-s)\alpha_{m-(s-j)}\beta_{n-j} = 0 \text{ and } (4.7.3b) \sum_{j=0}^{s}-j\alpha_{m-(s-j)}\beta_{n-j} = 0.$$

If $1 \le j \le s-1$ and $j \not\equiv 0 \pmod{p}$, then $\beta_{n-j} = 0$ by the induction hypothesis. If $1 \le j \le s-1$ and $j \equiv 0 \pmod{p}$, then $s-j < s$ and $s-j \not\equiv 0$. Hence, in this case, $\alpha_{m-(s-j)} = 0$ by the induction hypothesis. Thus, (4.7.3a) implies $-s\alpha_{m-s}\beta_n = 0$, which implies $\alpha_{m-s} = 0$, since $s \not\equiv 0 \pmod{p}$ and $\beta_n \ne 0$. By the same reasoning, (4.7.3b) implies $-s\alpha_m\beta_{n-s} = 0$, which implies $\beta_{n-s} = 0$ since $s \not\equiv 0 \pmod{p}$ and $\alpha_m \ne 0$. Therefore, if $\partial(f_1, f_2)/\partial(x_1, x_2) \in k^*$, then $\alpha_{m-i} = \beta_{n-i} = 0$ for each integer $i \not\equiv 0 \pmod{p}$. However, this leads to the following contradiction.

Since $\alpha_{m-i} = \beta_{n-i} = 0$ for each integer $i \not\equiv 0 \pmod{p}$ and since $m \not\equiv 0 \pmod{p}$ and $n \not\equiv 0 \pmod{p}$, we conclude $\alpha_0 = \beta_0 = 0$, but this contradicts what we observed above (second paragraph of case 2). Thus, Jacobian pairs satisfying Case 2 with $x_1 u_{x_1}$ and $x_2 u_{x_2}$ linearly independent over $k$ do not exist.

The next Proposition summarizes what we discovered in Example 4.7.

**Proposition 4.8.** *Let $u \in k[x_1, x_2]$ be nonzero and homogeneous of positive degree. Let $f_1 = x_1 h_1(u), f_2 = x_2 h_2(u)$, where $h_1, h_2 \in k[x]$. Assume $h_1$ and $h_2$ have positive degree. If $\partial(f_1, f_2)/\partial(x_1, x_2) \in k^*$, then $u = x^a y^b (h(x, y))^p$, for some positive integers $a$ and $b$ with $ab \not\equiv 0 \pmod{p}$ and some homogeneous $h(x_1, x_2) \in k[x_1, x_2]$. Furthermore, if $u = x^a y^b (h(x, y))^p$ as just described, then $\partial(f_1, f_2)/\partial(x_1, x_2) \in k^*$ if and only if $h_1(u)h_2(u) + auh_1'(u)h_2(u) + buh_1(u)h_2'(u) \in k^*$. If such is the case, then $[k(x_1, x_2) : k(f_1, f_2)] \equiv 0 \pmod{p}$ and at least one of the $f_i$ has degree greater than or equal to $p$.*

**Remark 4.9.** In this paper we have sought $k$-endomorphisms $\sigma$ of $k[x_1, x_2]$ such that $(\sigma(x_1), \sigma(x_2))$ is a Jacobian pair in order to test them against various characteristic $p > 0$ versions of the Jacobian Conjecture. So far, we have found two, namely, (a) $\sigma$ that are $p$-morphisms, and (b) $\sigma$ defined by $\sigma(x_i) = x_i h_i(u)$ as described in Proposition 4.8. The two endomorphisms described here are distinct in that each produces Jacobian pairs that the other doesn't. In fact, they produce entirely different sets of Jacobian pairs, since the pairs produced by the endomorphisms in (a) always have one point at infinity modulo $p$ by Theorem 2.18 and those produced by the endomorphisms in (b) always have two points at infinity modulo $p$ by Proposition 4.8. Of course, we can describe a third class of $k$ −endomorphisms $\sigma$ that produce Jacobian pairs $(\sigma(x_1), \sigma(x_2))$, namely, those that are compositional products of a finite number of the maps described in (a) and (b). This class of $k$ −endomorphism of $k[x_1, x_2]$, which is generated by those of type (a) and (b) above, produce some Jacobian pairs that those in (a) and (b) do not. For example, if we let $\tau$ be defined by $\tau(x_1) = x_1 + x_2$, $\tau(x_2) = x_1 - x_2$, and $\rho$ be the $k$ −endomorphism of $k[x_1, x_2]$ defined by $\rho(x_1) = x_1(\alpha x_1^{a-1} x_2^{mp-a} + 1)$ and $\rho(x_2) = x_2 \left( \frac{a\alpha}{a-1} x_1^{a-1} x_2^{mp-a} + 1 \right)$, which we introduced in Example 4.1, then $\sigma = \rho \circ \tau$ is a product of a type (a) and a type (b) morphism and each $\sigma(x_i)$ will have three points at infinity modulo $p$. Hence, $\sigma$ will be neither type $(a)$ or (b).

Three interrelated questions immediately arise from the above discussion: (1) Does the Separable Jacobian Conjecture hold for the class of $k$ −endomorphism of $k[x_1, x_2]$ that are compositional products of a finite number of type (a) and (b) morphisms? (2) Do the Low Degree Jacobian Conjectures hold for these endomorphisms? (3) Are all $k$ −endomorphisms of $k[x_1, x_2]$ that produce Jacobian pairs compositional products of a finite number of type (a) and (b) endomorphisms? We begin by answering Question 1.

**Proposition 4.10.** *Let $\sigma$ be a $k$ −endomorphism of $k[x_1, x_2]$ that is a compositional product of a finite number of type (a) and (b) endomorphisms. If $p$ does not divide $\deg(\sigma)$, then $\sigma$ is an automorphism.*

**Proof.** Let $\sigma$ be as described in the statement of the proposition. If one of the endomorphisms in the compositional product that defines $\sigma$ is type (b), then by Corollary 4.6 and Proposition 4.8, $\deg(\sigma) \equiv 0 \pmod{p}$. Hence, $\sigma$ is a $p$ −morphism. By Corollary 2.10, $\sigma$ is an automorphism of $k[x_1, x_2]$.

**Proposition 4.11.** *Let $\rho$ be a $k$ −endomorphism of $k[x_1, x_2]$ that is a compositional product of a finite number of type (a) and (b) endomorphisms. If $\rho$ is not a $k$ −automorphism of $k[x_1, x_2]$, then $\max\{\deg(\rho(x_1)), \deg(\rho(x_2))\} \geq p$.*

**Proof.** Recall that a type (a) morphism (i.e. $p$ −morphism) is a product of morphisms of type 1, type 2, and type 3 (Remark 2.20). If a morphism $\mu$ of $k[x_1, x_2]$ is type 1 or type 2, then $\deg(\mu) = 0$. If $\mu$ is type 3 or type (b), then $p$ divides $\deg(\mu)$ by Proposition 2.9 and Proposition 4.8. Hence, since $\rho$ is not a $k$ −automorphism of $k[x_1, x_2]$, the compositional product that defines $\rho$ must have at least one factor that is type 3 or type (b) by Corollary 2.6 and Corollary 2.10. On the other hand, if two or more of the endomorphisms in the compositional product that defines $\rho$ are either type 3 or type (b), then $\deg(\rho) \geq p^2$ by Corollary 2.6 and Corollary 2.10. Then, $p^2 \leq \deg(\rho(x_1)) \cdot \deg(\rho(x_2))$ by Proposition 1.5, which implies $\max\{\deg(\rho(x_1)), \deg(\rho(x_2))\} \geq p$. Hence, we are left to consider the case when exactly one factor in the compositional product that defines $\rho$ is type 3 or type (b). If it is type 3, then $\rho$ is a $p$ −morphism and by Proposition 2.23, $\max\{\deg(\rho(x_1)), \deg(\rho(x_2))\} \geq p$. Thus, we are only left with the case that exactly one factor in the compositional product that defines $\rho$ is type (b) and all others are type 1 and type 2. Then $\rho = \sigma_1 \phi \sigma_2$ where $\phi$ is a type (b) endomorphism and $\sigma_1$ and $\sigma_2$ are products of type 1 and type 2 maps and are therefore automorphisms.

We will proceed by induction on the number of type 1 and type 2 non-identity automorphisms in the compositional product that defines $\sigma_2$. If there are no such automorphisms (i.e. $\sigma_2$ is the identity map) then $\rho = \sigma_1 \phi$.

We have, $\sigma_1(x_i) = g_i(x_1, x_2)$, for some $g_i(x_1, x_2) \in k[x_1, x_2]$ for $i = 1, 2$.

Since $\phi$ is type (b), we have for each $i = 1, 2$, $\phi(x_i) = x_i h_i(u)$, for some $h_i \in k[x]$ of positive degree and $u = x_1^a x_2^b (h(x_1, x_2))^p$, with $a$ and $b$ positive integers such that $ab \not\equiv 0 \pmod{p}$ and for some $h(x_1, x_2)$ homogeneous in $k[x_1, x_2]$. We also have $a \cdot \deg(h_1) + b \cdot \deg(h_2) + 1 \equiv 0 \pmod{p}$ as we saw in Example 4.7. Without loss of generality, we may assume $\deg(h_1) \geq \deg(h_2)$. Then $\rho(x_1) = \sigma_1(\phi(x_1)) = \sigma_1\left(x_1 h_1(x_1^a x_2^b (h(x_1, x_2))^p)\right) = g_1 \cdot h_1(g_1^a g_2^b (h(g_1, g_2))^p)$. Thus, $\deg(\rho(x_1)) \geq \deg\left(g_1 \cdot h_1(g_1^a g_2^b)\right) = \deg(g_1) + \deg(h_1)(a \cdot \deg(g_1) + b \cdot \deg(g_2)) \geq \deg(g_1) + \deg(h_1)(a \cdot \deg(g_1) + b \cdot \deg(g_2)) = \deg(g_1) + a \cdot \deg(h_1) \cdot \deg(g_1) + b \cdot \deg(h_1) \cdot \deg(g_2) \geq \deg(g_1) + a \cdot \deg(h_1) \cdot \deg(g_1) + b \cdot \deg(h_2) \cdot \deg(g_2) \geq a \cdot \deg(h_1) + b \cdot \deg(h_2) + 1$, since $\deg(g_i) \geq 1$ for each $i$. Since $a \cdot \deg(h_1) + b \cdot \deg(h_2) + 1 \equiv 0 \pmod{p}$, we obtain, $\deg(\rho(x_1)) \geq p$. Hence, $\max\{\deg(\rho(x_1)), \deg(\rho(x_2))\} \geq p$.

From this point forward the proof will proceed exactly as in the last paragraph of Proposition 2.23. Hence, if $\rho$ is a $k$−endomorphism of $k[x_1, x_2]$ that is a compositional product of a finite number of type (a) and (b) endomorphisms and $\rho$ is not a $k$−automorphism of $k[x_1, x_2]$, then $\max\{\deg(\rho(x_1)), \deg(\rho(x_2))\} \geq p$.

**Remark 4.12.** From Proposition 4.10 and Proposition 4.11 we have affirmative answers to the first two questions listed in the closing paragraph of Remark 4.9. If the answer to the third question in Remark 4.9 is likewise affirmative, then the two variable Jacobian Conjecture holds in characteristic 0. Thus we anticipate that we will spend considerable time and effort in the future attempting to answer Question 3 of Remark 4.9. If we discover more examples of $k$−endomorphisms of $k[x_1, x_2]$ that produce Jacobian pairs yet are not finite products of Type (a) and Type (b) morphisms, we can revise Question (3) to include these morphisms in the potential generating set of morphisms producing Jacobian pairs.

---

[1] J. Lang, *A Jacobian identity in positive characteristic*, Journal of Commutative Algebra, Volume 7, Number 3, Fall 2015.

[2] O. H. Keller, *Ganze Cremona-Transformationen*, Montash Math. 47 (1939), 299-306.

[3] J. Lang and S. Mandal, *On Jacobian n-tuples in characteristic p*, Rocky Mountain Journal of Mathematics, Volume 23, Number 1, 1993.

[4] S. S. Abhyankar, *Expansion techniques in algebraic geometry*, Tata Lecture Notes, 1977.

[5] K. Adjamagboo, "On separable algebras over a U.F.D. and the Jacobian conjecture in any characteristic", *Automorphisms of affine spaces (Curaçao, 1994)*, Dordrecht: Kluwer Acad. Publ., pp. 89–103

[6] N. Jacobson, Basic Algebra II, Freeman & Company, W. H. (1980).

[7] H. Bass, E. Connell, and L. van den Dries, *The Jacobian conjecture: Reduction of degree and formal expansion of the inverse*, Bull. Amer. Math Soc. **7** (1982) 287-330.

[8] J. Lang and S. Mandal, *On Jacobian n-tuples in characteristic p*, Rocky Mountain Journal of Mathematics, Volume 23, Number 1, 1993.

[9] J. Lang, *A Jacobian identity in positive characteristic*, Journal of Commutative Algebra, Volume 7, Number 3, Fall 2015.
[10] H. W. E. Jung, *Uber ganze birationale Transformationen der Ebene*, J. Reine Angew. Math. 184 (1942) 161-174.
[11] W. Van der Kulk, On polynomial rings in two variables, Nieuw Arch. Wisk (3) (1953) 33-41.
[12] J. Lang, *Newton polygons of Jacobian pairs*, Journal of Pure and Applied algebra 72 (1991) 39-51.